\documentclass[AMA,STIX1COL]{WileyNJD-v2}
\usepackage{dsfont}
\usepackage{graphicx}
\usepackage{epstopdf}
\usepackage[caption=false]{subfig}
\usepackage{xcolor}
\usepackage{bm}
\newtheorem{Theorem}{Theorem}
\newtheorem{Lemma}{Lemma}

\newtheorem{Example}{Example}
\newtheorem{Remark}{Remark}

\newtheorem{Assumption}{Assumption}
\allowdisplaybreaks
\articletype{Research Article}%
\begin{document}

\title{Composite Optimization with Coupling Constraints via Dual Proximal Gradient Method with Applications to Asynchronous Networks
}

\author[1]{Jianzheng Wang}

\author[1]{Guoqiang Hu}

\authormark{Jianzheng Wang \textsc{et al}}

\address[]{\orgdiv{School of Electrical and Electronic Engineering}, \orgname{Nanyang Technological University}, \orgaddress{\state{639798}, \country{Singapore}}}

\corres{Guoqiang Hu, School of Electrical and Electronic Engineering, Nanyang Technological University, 639798, Singapore. \email{gqhu@ntu.edu.sg}}


\abstract[Summary]{In this paper, we consider solving a composite optimization problem with affine coupling constraints in a multi-agent network based on proximal gradient method. In this problem, all the agents jointly minimize the sum of individual cost functions composed of smooth and possibly non-smooth parts. To this end, we derive the dual problem by the concept of Fenchel conjugate, which gives rise to the dual proximal gradient algorithm by allowing for the asymmetric individual interpretations of the global constraints. Then, an asynchronous dual proximal gradient algorithm is proposed for the asynchronous networks with heterogenous step-sizes and communication delays. For both the two algorithms, if the non-smooth parts of the objective functions are simple-structured, we only need to update dual variables by some simple operations, accounting for the reduction of the overall computational complexity. Analytical convergence rate of the proposed algorithms is derived and their efficacy is verified by solving a social welfare optimization problem of electricity market in the numerical simulation.}


\maketitle

\section{Introduction}\label{}

\subsection{Background and Motivation}

{D}{ecentralized} optimization has drawn much attention due to its prominent advantage in solving various mathematical optimization problems with large data set and decentralized decision variables in multi-agent networks.\cite{luo2014provably,lee2017speeding,bai2017distributed} In those problems, each agent usually maintains a local decision variable, and the optimal solution of the system is achieved through multiple rounds of communications and strategy-makings.\cite{yuan2013distributed} In this work, we consider a class of optimization problems with composite cost functions, i.e., composed of smooth (differentiable) and possibly non-smooth (non-differentiable) parts, arising from various fields, such as Lasso regressions, resource allocation problems and support vector machines. \cite{hans2009bayesian,beck2014fast,zhao2017scope}

To solve those problems, most existing works require the update of primal variables with some costly computations, which increase the overall computational complexity. Meanwhile, with the presence of the asynchrony of large-scale networks in various fields, more explorations on asynchronous optimization algorithms are needed.\cite{bertsekas1989parallel} As widely discussed, proximal gradient based algorithms can take the advantage of some simple-structured cost functions and are usually numerically more stable than the subgradient based counterparts.\cite{bertsekas2011incremental} With the above motivation, in this work, we aim to develop an efficient optimization algorithm for decentralized optimization problems (DOPs) based on proximal gradient method and further investigate its efficacy in asynchronous networks.

\subsection{Literature Review}

In this work, we focus on optimizing a class of composite DOPs subject to affine coupling constraints. To solve these problems, applicable techniques include primal-dual subgradient methods,\cite{nesterov2009primal} alternating direction method of multipliers,\cite{wang2019global} and proximal gradient methods,\cite{parikh2014proximal} etc. DOPs with coupling constraints are actively investigated in the recent works,\cite{mateos2016distributed,necoara2017fully,notarnicola2019constraint,li2020distributed,falsone2016distributed,falsone2017dual,zhu2011distributed,chang2016proximal,chang2014distributed,simonetto2016primal} where the optimal solution to the primal problems is usually achieved with the update of both primal and dual variables.
An alternative solution, as discussed by Notarnicola et al.,\cite{notarnicola2017duality,notarnicola2018duality} is resorting to the dual problems, where the computation on the primal variables is not required. However, the algorithms in References 23 and 24 involve some inner-loop optimization processes, which increase the overall computational complexity if the primal cost functions possess some non-smooth characteristics. To further improve the computational efficiency, dual proximal gradient (DPG) methods for solving composite optimization problems were investigated recently,\cite{notarnicola2016asynchronous,beck2014fast,kim2016fast} where, however, no general affine coupling constraint was considered.

To explore some efficient decentralized algorithms, different from the existing works, the new features of this work are twofold. First, to the best knowledge of the authors, this is the first work that investigates DPG method with general affine constraints with specific network typologies. By the proposed DPG algorithm, the updating of the primal variables is not compulsory. Furthermore, if the proximal mapping of the non-smooth parts in the primal problem can be explicitly given, we only need to update the dual variables by some simple operations,\footnote[1]{We also note that some dual algorithms dealing with smooth cost functions also can avoid the update of primal variables.\cite{necoara2015linear} However, directly extending their algorithms to non-smooth cases can be costly in the sense that the computation of the gradient of the conjugate of a non-smooth function can be costly in general. Therefore, the contribution to the computational efficiency (also Asyn-DPG algorithm as introduced later) is established for possibly non-smooth cost functions.} e.g., basic proximal mappings and gradient based iterations, which technically can be more efficient than the existing algorithms with some costly computations on the primal variables or other auxiliary variables.\cite{mateos2016distributed,necoara2017fully,notarnicola2019constraint,li2020distributed,falsone2016distributed,falsone2017dual,zhu2011distributed,chang2016proximal,chang2014distributed,simonetto2016primal} As another feature, the asymmetric individual interpretation of the agents on the global constraints is considered, where no uniform knowledge of the global constraints is required.

Second, we propose an asynchronous dual proximal gradient (Asyn-DPG) algorithm, which can be viewed as an extension of DPG algorithm by considering heterogenous step-sizes and communication delays. Specifically, the outdated information is addressed through deterministic analysis,\cite{zhou2018distributed,chang2016asynchronous,hale2017asynchronous,nedic2001incremental,cannelli2020asynchronous,tian2020achieving,li2013parameter,hong2017distributed} which is advantageous over some stochastic models \cite{cannelli2019asynchronous} in the sense that the probability distribution of random factors can be difficult to acquire in some problems and may introduce inaccuracy issues due to limited historical data.\cite{bertsekas2002introduction}
However, the problem setup in References 28, 29, 31-35 either only considers certain special form of affine coupling constraints or does not incorporate any coupling constraint. In addition, the algorithms discussed in References 30 and 33 dealing with smooth cost functions will hamper their usage in non-smooth optimization problems. Different from all the aforementioned works, we will show that if the upper bound of communication delays is finite and the non-smooth parts in the primal problem are simple-structured, we only need to update dual variables with some simple operations, which is still a distinct advantage to reduce the computational complexity.

We hereby summarize the contributions of this work as follows.
\begin{itemize}
  \item We consider a class of composite DOPs with both local convex and affine coupling constraints. To solve these problems, a DPG algorithm is proposed by formulating the dual problems. Then, an Asyn-DPG algorithm is built upon the structure of DPG algorithm, which can be applied to asynchronous networks with heterogenous step-sizes and communication delays. In addition, the asymmetric individual interpretations of the global constraints are considered, which is more adaptive to large-scale networks in the sense that no uniform knowledge of the global constraints for the agents is required.
  \item Provided that the non-smooth parts of the cost functions in the primal DOPs are with some simple structures, the proposed DPG and Asyn-DPG algorithms only require the update of dual variables with some simple operations, accounting for the reduction of the overall computational complexity. In addition, our algorithms require some commonly used assumptions on the primal problems and explicit convergence rates are provided for all the discussed scenarios.
\end{itemize}

\subsection{Paper Structure and Notations}

The remainder of this paper is organized as follows. Section \ref{sa2} presents some frequently used definitions in this work and their properties. Section \ref{sa3} formulates the primal problem of interest and gives some basic assumptions. In Section \ref{sa4}, two proximal gradient algorithms, namely DPG and Asyn-DPG, are proposed based on different network settings. The convergence analysis of the discussed algorithms is conducted in Section \ref{sa5}. The efficacy of the proposed algorithms is verified by a numerical simulation in Section \ref{sa6}. Section \ref{sa7} concludes this paper.

$\mathds{N}$ and $\mathds{N}_+$ denote the non-negative and positive integer spaces, respectively. Let notation $\mid \mathcal{S}\mid$ be the size of set $\mathcal{S}$.
Operator $(\cdot)^{\top}$ represents the transpose of a matrix.
$\mathcal{S}_1 \times \mathcal{S}_2$ denotes the Cartesian product of sets $\mathcal{S}_1$ and $\mathcal{S}_2$.
$\mathbf{relint}\mathcal{S}$ represents the relative interior of set $\mathcal{S}$. Let $\lfloor u \rfloor$ ($\lceil u \rceil $) be the largest integer smaller than (smallest integer no smaller than) scalar $u$. $\| \cdot\|_1$ and $\| \cdot \|$ refer to the $l_1$ and $l_2$-norms, respectively. $\langle \cdot,\cdot \rangle$ is an inner product operator. $\otimes$ is the Kronecker product operator. $\mathbf{0}_n$ and $\mathbf{1}_n$ denote the $n$-dimensional column vectors with all elements being 0 and 1, respectively. $\mathbf{I}_n$ and $\mathbf{O}_{n \times m}$ denote the $n$-dimensional identity matrix and $(n \times m)$-dimensional zero matrix, respectively.

\section{Preliminaries}\label{sa2}

In this section, we present some fundamental definitions and properties of graph theory, proximal mapping, and Fenchel conjugate.

\subsection{Graph Theory}\label{d2}
A multi-agent network can be described by an undirected graph ${\mathcal{G}}= \{{\mathcal{V}},{\mathcal{E}}\}$, which is composed of the set of vertices ${\mathcal{V}} = \{1,2,...,N\}$ and set of edges ${\mathcal{E}} \subseteq \{ (i,j)| i,j \in \mathcal{V} \hbox{ and } i \neq j\}$ with $(i,j) \in \mathcal{E}$ an unordered pair (no self-loop). A graph $\mathcal{G}$ is said connected if there exists at least one path between any two distinct vertices. A graph $\mathcal{G}$ is said fully connected if there is a unique edge between any two distinct vertices. ${\mathcal{V}}_i =  \{ j | (i,j) \in \mathcal{E}\}$ denotes the set of the neighbours of agent $i$. Let ${\mathbf{L}} \in \mathds{R}^{N \times N}$ denote the Laplacian matrix of ${\mathcal{G}}$. Let $d_{ij}$ be the element at the cross of the $i$th row and $j$th column of ${\mathbf{L}}$. Thus, $d_{ij} = -1$ if $(i,j) \in {\mathcal{E}}$, $d_{ii} = \mid {\mathcal{V}}_i \mid$, and $d_{ij} = 0$ otherwise, $i,j \in {\mathcal{V}}$.\cite{chung1997spectral}

\subsection{Proximal Mapping}
A proximal mapping of a proper, closed, and convex function $\zeta: \mathds{R}^n \rightarrow (-\infty,+\infty]$ is defined by $\mathrm{prox}^{\alpha}_{\zeta} [\mathbf{u}] =  \arg \min_{\mathbf{v}} ( \zeta(\mathbf{v}) + \frac{1}{2{\alpha}} \| \mathbf{v} - \mathbf{u}\|^2 )$ with ${\alpha}>0$ and $\mathbf{u} \in \mathds{R}^n$.\cite{notarnicola2016asynchronous} \footnote[2]{The proximal mapping can be equivalently written as $\mathrm{prox}_{{\alpha} \zeta}$ as in some other works. Hence, we will not mention this difference for conciseness when citing those works in the rest of this paper.}

\subsection{Fenchel Conjugate}\label{de1}
Let $\zeta: \mathds{R}^n \rightarrow (-\infty,+\infty]$ be a proper function. The Fenchel conjugate of $\zeta$ is defined as $\zeta^{\diamond}(\mathbf{u})= \sup_{\mathbf{v}} \{\mathbf{u}^{\top} \mathbf{v}-\zeta(\mathbf{v})\}$, which is convex.\cite[Sec. 3.3]{borwein2010convex}

\begin{Lemma}\label{md}
(Extended Moreau Decomposition\cite[Thm. 6.45]{beck2017first}) Let $\zeta: \mathds{R}^n \rightarrow (-\infty,+\infty]$ be a proper, closed, convex function and $\zeta^{\diamond}$ be its Fenchel conjugate. Then, for all $\mathbf{u} \in \mathds{R}^n$ and $\alpha >0$, we have
\begin{align}\label{am}
 \mathbf{u} = \alpha \mathrm{prox}^{\frac{1}{\alpha}}_{\zeta^{\diamond}} [\frac{\mathbf{u}}{\alpha}] + \mathrm{prox}^{\alpha}_{\zeta} [\mathbf{u}].
\end{align}
\end{Lemma}

\begin{Lemma}\label{l1}
Let $\zeta: \mathds{R}^n \rightarrow (-\infty,+\infty]$ be a proper, closed, $\sigma$-strongly convex function and $\zeta^{\diamond}$ be its Fenchel conjugate, $\sigma >0$. Then,
\begin{align}
\nabla_{\mathbf{u}} \zeta^{\diamond}(\mathbf{u}) = \arg \max\limits_{\mathbf{v}} ( \mathbf{u}^{\top} \mathbf{v} - \zeta(\mathbf{v})),
\end{align}
and $\nabla \zeta^{\diamond}$ is Lipschitz continuous with constant $\frac{1}{\sigma}$.\cite[Lemma V.7]{notarnicola2016asynchronous}
\end{Lemma}

\section{Problem Formulation}\label{sa3}

The considered optimization problem and relevant assumptions are presented in this section.

Consider a multi-agent network $\mathcal{G}= \{\mathcal{V},\mathcal{E}\}$ and a global cost function $H(\mathbf{x})= \sum_{i \in \mathcal{V}} H_i(\mathbf{x}_i)$, $\mathbf{x}_i \in \mathds{R}^M$, $\mathbf{x}=[\mathbf{x}^{\top}_1,...,\mathbf{x}^{\top}_{N}]^{\top} \in \mathds{R}^{{N}M}$. Agent $i$ maintains a private cost function $H_i(\mathbf{x}_i)= f_i(\mathbf{x}_i)+g_i(\mathbf{x}_i)$. Let $\Omega_i \subseteq \mathds{R}^M$ be the feasible region of $\mathbf{x}_i$. Then, the feasible region of $\mathbf{x}$ can be defined by $\Omega=  \Omega_1 \times \Omega_2 \times ... \times \Omega_{N} \subseteq \mathds{R}^{NM}$. We consider a global affine constraint $\mathbf{A}\mathbf{x} = \mathbf{b}$, $\mathbf{A} \in \mathds{R}^{B \times {N}M}$, $\mathbf{b} \in \mathds{R}^{B}$. Then, a DOP of $\mathcal{V}$ can be formulated as
\begin{align}\label{}
\textbf{(P1)} \quad  \min\limits_{\mathbf{x} \in \Omega} \quad  & \sum_{i \in \mathcal{V}} H_i(\mathbf{x}_i) \nonumber \\
 \hbox{subject to} \quad & \mathbf{A}\mathbf{x} = \mathbf{b}. \nonumber
\end{align}

\begin{Assumption}\label{a0}
$\mathcal{G}$ is undirected and connected.
\end{Assumption}

\begin{Assumption}\label{a1}
$f_i: \mathds{R}^M \rightarrow (-\infty,+\infty]$ is a proper, closed, differentiable, and $\sigma_i$-strongly convex extended real-valued function, $\sigma_i > 0$; $g_i: \mathds{R}^M \rightarrow (-\infty,+\infty]$ is a proper, closed and convex extended real-valued function, $i\in \mathcal{V}$.
\end{Assumption}

The assumptions in Assumption \ref{a1} are often discussed in composite optimization problems.\cite{shi2015proximal,wang2020distributed,schmidt2011convergence,chang2014multi,florea2020generalized,beck2014fast,notarnicola2016asynchronous}

\begin{Assumption}\label{a1+1}
(Constraint Qualification) $\Omega_i$ is non-empty, convex and closed, $i\in \mathcal{V}$; there exists an $\breve{\mathbf{x}} \in \mathbf{relint}\Omega$ such that $\mathbf{A}\breve{\mathbf{x}} = \mathbf{b}$.\cite{boyd2004convex}
\end{Assumption}

In the following, we consider that each agent $i$ maintains a private constraint $\mathcal{X}_i= \{\mathbf{x} \in \mathds{R}^{{N}M}|\mathbf{A}^{(i)}\mathbf{x} = \mathbf{b}^{(i)}\}$, which can be regarded as an individual interpretation of the global constraint $\mathcal{X}= \{\mathbf{x} \in \mathds{R}^{{N}M}| \mathbf{A}\mathbf{x} = \mathbf{b}\}$, $\mathbf{A}^{(i)} \in \mathds{R}^{B \times {N}M}$, $\mathbf{b}^{(i)} \in \mathds{R}^{B}$. Therefore, it is reasonable to assume that $\bigcap_{i\in \mathcal{V}} \mathcal{X}_i = \mathcal{X}$. Then, Problem (P1) can be equivalently written as
\begin{align}\label{}
\textbf{(P2)} \quad  \min\limits_{\mathbf{x}} \quad  & \sum_{i \in \mathcal{V}} (H_i(\mathbf{x}_i) +\mathds{I}_{\Omega_{i}}(\mathbf{x}_i)) \nonumber \\
 \hbox{subject to} \quad & \mathbf{A}^{(i)}\mathbf{x} = \mathbf{b}^{(i)},  \quad  \forall i \in \mathcal{V}, \nonumber
\end{align}
with $\mathds{I}_{\Omega_{i}} (\mathbf{x}_i)= \left\{\begin{array}{ll}
                    0, & \hbox{if $\mathbf{x}_i \in \Omega_i$}, \\
                    +\infty,  & \hbox{otherwise}
                  \end{array}
                  \right.$.\cite{boyd2004convex}

To facilitate the following discussion, we let $\mathbf{A}^{(l)}_i \in \mathds{R}^{B \times M}$ denote the $i$th column sub-block of $\mathbf{A}^{(l)}$, i.e., $\mathbf{A}^{(l)} = [\mathbf{A}^{(l)}_1,...,\mathbf{A}^{(l)}_i,...,\mathbf{A}^{(l)}_{N}]$, $i,l \in \mathcal{V}$.

\begin{Assumption}\label{as4}
Assume that $\mathbf{A}^{(i)}\mathbf{x} = \mathbf{b}^{(i)}$ only contains the decision variables of agent $i$ and its neighbours, i.e.,
\begin{align}\label{37}
\mathbf{A}^{(i)}_j=\mathbf{O}_{B \times M}, \quad \forall (i,j) \notin \mathcal{E} \hbox{ and } i \neq j.
\end{align}
\end{Assumption}
Some DOPs complying with Assumption \ref{as4} will be discussed in Section \ref{6}.


\section{Dual Proximal Gradient Based Algorithm Development}\label{sa4}

In this section, we will develop two dual proximal gradient based algorithms for solving the problem of interest under different assumptions on networks.

\subsection{Dual Problem}

By introducing a slack vector $\mathbf{z}=[\mathbf{z}^{\top}_1,...,\mathbf{z}^{\top}_{N}]^{\top}$, Problem (P2) can be equivalently written as
\begin{align}\label{}
 \textbf{(P3)} \quad    \min\limits_{\mathbf{x},\mathbf{z}} \quad  & \sum_{i \in \mathcal{V}} (f_i(\mathbf{x}_i) + (g_i+ \mathds{I}_{\Omega_{i}}) (\mathbf{z}_i)) \nonumber \\
  \hbox{subject to} \quad  & \mathbf{A}^{(i)}\mathbf{x} = \mathbf{b}^{(i)}, \quad  \mathbf{x}_i = \mathbf{z}_i, \quad \forall i \in \mathcal{V}. \nonumber
\end{align}
Then, the Lagrangian function of Problem (P3) can be written as
\begin{align}\label{28}
L (\mathbf{x}, \mathbf{z}, \bm{\theta},\bm{\mu})
= & \sum_{i \in \mathcal{V}} (f_i(\mathbf{x}_i) + (g_i+{\mathds{I}}_{{\Omega}_{i}}) (\mathbf{z}_i) + \bm{\mu}_i^{\top} (\mathbf{x}_i - \mathbf{z}_i) ) + \sum_{i \in \mathcal{V}} \bm{\theta}_i^{\top}(\mathbf{A}^{(i)}\mathbf{x} - \mathbf{b}^{(i)} ) \nonumber \\
= & \sum_{i \in \mathcal{V}} (f_i(\mathbf{x}_i) + (g_i+{\mathds{I}}_{{\Omega}_{i}}) (\mathbf{z}_i) + \bm{\mu}_i^{\top} (\mathbf{x}_i - \mathbf{z}_i) )   + \sum_{i \in \mathcal{V}} \mathbf{x}^{\top}_i \sum_{l\in \mathcal{V}}  (\mathbf{A}_i^{(l)})^{\top} \bm{\theta}_l  -  \sum_{i \in \mathcal{V}} (\mathbf{b}^{(i)} )^{\top}\bm{\theta}_i \nonumber \\
= & \sum_{i \in \mathcal{V}} (f_i(\mathbf{x}_i) +  \mathbf{x}_i^{\top} ( \sum_{l\in \mathcal{V}} (\mathbf{A}^{(l)}_i )^{\top} \bm{\theta}_l + \bm{\mu}_i )   + (g_i+{\mathds{I}}_{{\Omega}_{i}}) (\mathbf{z}_i)- \mathbf{z}_i^{\top} \bm{\mu}_i - (\mathbf{b}^{(i)} )^{\top}\bm{\theta}_i) ,
\end{align}
where we use
\begin{align}\label{}
\sum_{i \in \mathcal{V}} \bm{\theta}_i^{\top} \mathbf{A}^{(i)}\mathbf{x}
& =  \sum_{i \in \mathcal{V}} \sum_{l\in \mathcal{V}} \bm{\theta}_i^{\top} \mathbf{A}_l^{(i)}\mathbf{x}_l =  \sum_{i \in \mathcal{V}} \sum_{l\in \mathcal{V}} \bm{\theta}_l^{\top} \mathbf{A}_i^{(l)}\mathbf{x}_i   \nonumber \\
& =  \sum_{i \in \mathcal{V}} \sum_{l\in \mathcal{V}} \mathbf{x}^{\top}_i  (\mathbf{A}_i^{(l)})^{\top} \bm{\theta}_l   = \sum_{i \in \mathcal{V}} \mathbf{x}^{\top}_i \sum_{l\in \mathcal{V}}  (\mathbf{A}_i^{(l)})^{\top} \bm{\theta}_l
\end{align}
with ${\bm{\theta}} = [\bm{\theta}_1^{\top},...,\bm{\theta}_N^{\top}]^{\top} \in \mathds{R}^{{N}B}$ and $ {\bm{\mu}}=[\bm{\mu}_1^{\top},...,\bm{\mu}_N^{\top}]^{\top} \in \mathds{R}^{{N}M}$. $\bm{\theta}_i$ and $\bm{\mu}_i$ denote the Lagrangian multiplier vectors associated with constraints $\mathbf{A}^{(i)}\mathbf{x} = \mathbf{b}^{(i)}$ and $\mathbf{x}_i = \mathbf{z}_i$, respectively.

Therefore, the dual function can be obtained by minimizing ${L}(\mathbf{x},\mathbf{z}, {\bm{\theta}},{\bm{\mu}})$ with $(\mathbf{x},\mathbf{z})$, which gives
\begin{align}\label{v1}
 {V}  ({\bm{\theta}},  {\bm{\mu}}) =  & \min\limits_{{\mathbf{x}},{\mathbf{z}}} \sum_{i \in \mathcal{V}} (f_i(\mathbf{x}_i) +  \mathbf{x}_i^{\top} ( \sum_{l\in \mathcal{V}} (\mathbf{A}^{(l)}_i )^{\top} \bm{\theta}_l + \bm{\mu}_i ) + (g_i+{\mathds{I}}_{{\Omega}_{i}}) (\mathbf{z}_i)- \mathbf{z}_i^{\top} \bm{\mu}_i - (\mathbf{b}^{(i)} )^{\top}\bm{\theta}_i )\nonumber \\
 = &  \min\limits_{{\mathbf{x}},{\mathbf{z}}} \sum_{i \in \mathcal{V}} ( f_i(\mathbf{x}_i) -  \mathbf{x}_i^{\top} {\mathbf{C}}_i {\bm{\lambda}} + (g_i + {\mathds{I}}_{{\Omega}_{i}}) (\mathbf{z}_i)- \mathbf{z}_i^{\top} {\mathbf{F}}_i  {\bm{\lambda}} - {\mathbf{E}}_i {\bm{\lambda}} )  \nonumber \\
 = & \sum_{i \in \mathcal{V}} (- f^{\diamond}_i( {\mathbf{C}}_i {\bm{\lambda}} ) - {\mathbf{E}}_i {\bm{\lambda}} - ( g_i +{\mathds{I}}_{{\Omega}_{i}} )^{\diamond}({\mathbf{F}}_i  {\bm{\lambda}})) ,
\end{align}
where
\begin{align}\label{}
  {\mathbf{C}}_i & =  [ -(\mathbf{A}^{(1)}_i)^{\top} ,..., -(\mathbf{A}^{({N})}_i)^{\top},  \mathbf{O}_{M \times (i-1)M}, -\mathbf{I}_{M}, \mathbf{O}_{M \times ({N}-i)M}] \in \mathds{R}^{M \times ({N}M+{N}B)}, \\
  {\mathbf{F}}_i & =  [\mathbf{O}_{M \times ({N}B + (i-1)M)}, \mathbf{I}_{M}, \mathbf{O}_{M \times ({N}-i)M}] \in \mathds{R}^{M\times ({N}M+{N}B)}, \\
  {\mathbf{E}}_i & = [\mathbf{0}^{\top}_{(i-1)B}, (\mathbf{b}^{(i)})^{\top}, \mathbf{0}^{\top}_{(N-i)B+NM}] \in \mathds{R}^{1 \times ({N}M+{N}B)}, \\
  {\bm{\lambda}} & = [{\bm{\theta}}^{\top}, {\bm{\mu}}^{\top}]^{\top} \in \mathds{R}^{{N}M+{N}B}.
\end{align}
Then, the dual problem of Problem (P3) can be formulated as
\begin{align}
\textbf{(P4)} \quad \min\limits_{{\bm{\lambda}}} \quad {\Psi}({\bm{\lambda}})  \nonumber
\end{align}
where
\begin{align}\label{}
  {\Psi}({\bm{\lambda}}) & = {P}({\bm{\lambda}}) + {Q}({\bm{\lambda}}), \\
  {P} ({\bm{\lambda}}) & = \sum_{i\in \mathcal{V}} {p}_i ({\bm{\lambda}}) , \\
  {Q} ({\bm{\lambda}}) &  = \sum_{i\in \mathcal{V}} {q}_i ({\bm{\mu}_i}), \\
  {p}_i ({\bm{\lambda}}) & = f^{\diamond}_i( {\mathbf{C}}_i {\bm{\lambda}} )  + {\mathbf{E}}_i {\bm{\lambda}}, \\
{q}_i ({\bm{\mu}_i}) &= (g_i +{\mathds{I}}_{{\Omega}_{i}} )^{\diamond}({\mathbf{F}}_i{\bm{\lambda}}) = (g_i +{\mathds{I}}_{{\Omega}_{i}} )^{\diamond}({\bm{\mu}_i}).
\end{align}
Define $\mathcal{H}$ as the set of the optimal solutions to Problem (P4).

\subsection{Discussion on Constraints}\label{6}

Assumption \ref{as4} can be satisfied in some DOPs. See Examples \ref{ex0} to \ref{ex3}.

\begin{Example}\label{ex0}
$\mathcal{G}$ is fully connected.
\end{Example}

\begin{Example}\label{ex1}
In some applications, e.g., telecommunication and machine learning problems, $\mathbf{A}^{(i)}\mathbf{x} = \mathbf{b}^{(i)}$ can be defined by an edge-constraint maintained by agent $i$.\cite{zhang2017distributed}
\end{Example}

\begin{Example}\label{ex2}
If $\mathbf{A}_i^{(i)}=|\mathcal{V}_i| \mathbf{I}_{M}, \mathbf{A}_l^{(i)} =-\mathbf{I}_{M}$, $\mathbf{A}_{j}^{(i)}=\mathbf{O}_{M \times M}$, $\mathbf{b}^{(i)}= \mathbf{0}$, $\forall (i,l) \in \mathcal{E}$, $(i,j) \notin \mathcal{E} \hbox{ and } i \neq j$, then the constraints of Problem (P3) can be written as $\bm{\mathsf{L}} \mathbf{x} = \mathbf{0}$ with $\bm{\mathsf{L}} = {\mathbf{L}} \otimes \mathbf{I}_{M} \in \mathds{R}^{N M \times  N M}$, which means Problem (P3) essentially is a consensus optimization problem.\cite{sun2017distributed,wen2012consensus,yuan2015gradient}
\end{Example}

\begin{Example}\label{ex3}
Consider a set of consensus constraints of agent $i$: $\mathbf{x}_i - \mathbf{x}_j = \mathbf{0}$, $\forall j \in \mathcal{V}_i$.\cite{notarnicola2016asynchronous} Then, for any agent pair $(i,j) \in \mathcal{E}$, the individual constraints of agents $i$ and $j$ include $\mathbf{x}_i - \mathbf{x}_j = \mathbf{0}$ and $\mathbf{x}_j - \mathbf{x}_i = \mathbf{0}$, respectively. Therefore, the asymmetric constraints can be viewed as a generalization of the asymmetric consensus constraints discussed in this example.
\end{Example}

\begin{Remark}
In Examples 1-4, the asymmetric constraints are more adaptive to large-scale networks in the sense that establishing a global $\mathbf{A}\mathbf{x} = \mathbf{b}$ by integrating the overall decentralized or even distributed constraints may be costly, especially when the network sizes and individual constraints vary constantly.\cite{franceschelli2020stability,chen2018bandit} \footnote[3]{We assume the parameters of the networks are invariant when optimizations are in progress.} For example, when certain agent $i$ joins the network, he only needs to broadcast $\mathbf{A}^{(i)}$ to neighbours such that ${\bm{\lambda}}$ can be augmented directly as in Problem (P3), without changing the network-wide constraint architecture seriously by rebuilding $\mathbf{A}\mathbf{x} = \mathbf{b}$.
\end{Remark}

In practice, the asymmetric individual constraints can be generated by interpreting some common global constraints by user-defined linear transformations. For instance, agent $i$ may interpret constraint $\mathbf{A}\mathbf{x} = \mathbf{b}$ by transformation $\mathbf{T}_i\mathbf{A}\mathbf{x} = \mathbf{T}_i\mathbf{b}$, i.e., $\mathbf{A}^{(i)}=\mathbf{T}_i\mathbf{A}$ and $\mathbf{b}^{(i)}=\mathbf{T}_i\mathbf{b}$. See Example \ref{ex4}.

\begin{Example}\label{ex4}
Consider a global affine constraint $\mathcal{X} = \left\{ \mathbf{x} \in \mathds{R}^3 \bigg| \left[
  \begin{array}{ccc}
    1 & 1 & 0 \\
    2 & 0 & 1 \\
  \end{array}
\right]
\mathbf{x} =
\left[
  \begin{array}{c}
    1 \\
    2 \\
  \end{array}
\right] \right\}$ for a 3-agent network. The individual constraints maintained by agents 1, 2, and 3 are assumed to be
\begin{align}\label{}
   \mathcal{X}_1 = & \left\{ \mathbf{x} \in \mathds{R}^3 \bigg|\left[
  \begin{array}{ccc}
    -1 & -1 & 0 \\
    1 & 0 & \frac{1}{2} \\
  \end{array}
\right]
\mathbf{x} =
\left[
  \begin{array}{c}
    -1 \\
    1  \\
  \end{array}
\right] \right\}, \\
\mathcal{X}_2 = & \left\{ \mathbf{x} \in \mathds{R}^3 \big| \left[
  \begin{array}{ccc}
    2 & 2 & 0 \\
  \end{array}
\right]
\mathbf{x} = 2 \right\}, \\
\mathcal{X}_3 = & \left\{ \mathbf{x} \in \mathds{R}^3 \big| \left[
  \begin{array}{ccc}
    -2 & 0 & -1 \\
  \end{array}
\right]
\mathbf{x} = -2 \right\},
\end{align}
respectively, where $\mathcal{X} = \mathcal{X}_1  \bigcap \mathcal{X}_2  \bigcap \mathcal{X}_3 $. In this example, $\mathbf{T}_1=\left[
  \begin{array}{cc}
    -1 & 0 \\
    0 & \frac{1}{2}
  \end{array}
\right]$, $\mathbf{T}_2=\left[
  \begin{array}{cc}
    2 & 0
  \end{array}
\right]$, and $\mathbf{T}_3=\left[
  \begin{array}{cc}
    0 & -1
  \end{array}
\right]$.
\end{Example}

\subsection{Dual Proximal Gradient Algorithm}

In this subsection, we propose a DPG algorithm to solve Problem (P4).
The DPG algorithm is designed as
\begin{align}\label{8}
 {\bm{\lambda}}(k+1) = & \mathrm{prox}^{c}_{{Q}}[{\bm{\lambda}}(k) - c  \nabla_{{\bm{\lambda}}}  {P} ({\bm{\lambda}}(k))],
\end{align}
which means
\begin{align}\label{55-1}
 \left[
  \begin{array}{c}
    {\bm{\theta}_1}(k+1) \\
        \vdots \\
    {\bm{\theta}_N}(k+1)  \\
        \bm{\mu}_1(k+1) \\
    \vdots \\
    \bm{\mu}_N(k+1) \\
  \end{array}
\right]
= & \left[
            \begin{array}{c}
               {\bm{\theta}_1}(k) - c  \nabla_{{\bm{\theta}_1}} {P} ({\bm{\lambda}}(k)) \\
               \vdots \\
               {\bm{\theta}_N}(k) - c  \nabla_{ {\bm{\theta}_N}} {P} ({\bm{\lambda}}(k)) \\
               \mathrm{prox}^{c}_{{q}_1}[{\bm{\mu}_1}(k) - c  \nabla_{{\bm{\mu}_1}} {P} ({\bm{\lambda}}(k))] \\
                              \vdots \\
               \mathrm{prox}^{c}_{{q}_{N}}[{\bm{\mu}_N}(k) - c  \nabla_{{\bm{\mu}_N}} {P} ({\bm{\lambda}}(k))] \\
            \end{array}
          \right]
\end{align}
with $ \nabla_{{\bm{\lambda}}}  {P} = [\nabla^{\top}_{{\bm{\theta}_1}} {P},...,\nabla^{\top}_{{\bm{\theta}_N}} {P},\nabla^{\top}_{{\bm{\mu}_1}} {P} ,...,\nabla^{\top}_{{\bm{\mu}_N}} {P}]^{\top}$ and $c>0$, $k \in \mathds{N}$. The proximal mapping for computing ${\bm{\theta}}$ is omitted since ${\bm{\theta}}$ is not contained by ${Q}$.

\begin{figure}[htbp]
  \centering
  \includegraphics[width=7cm]{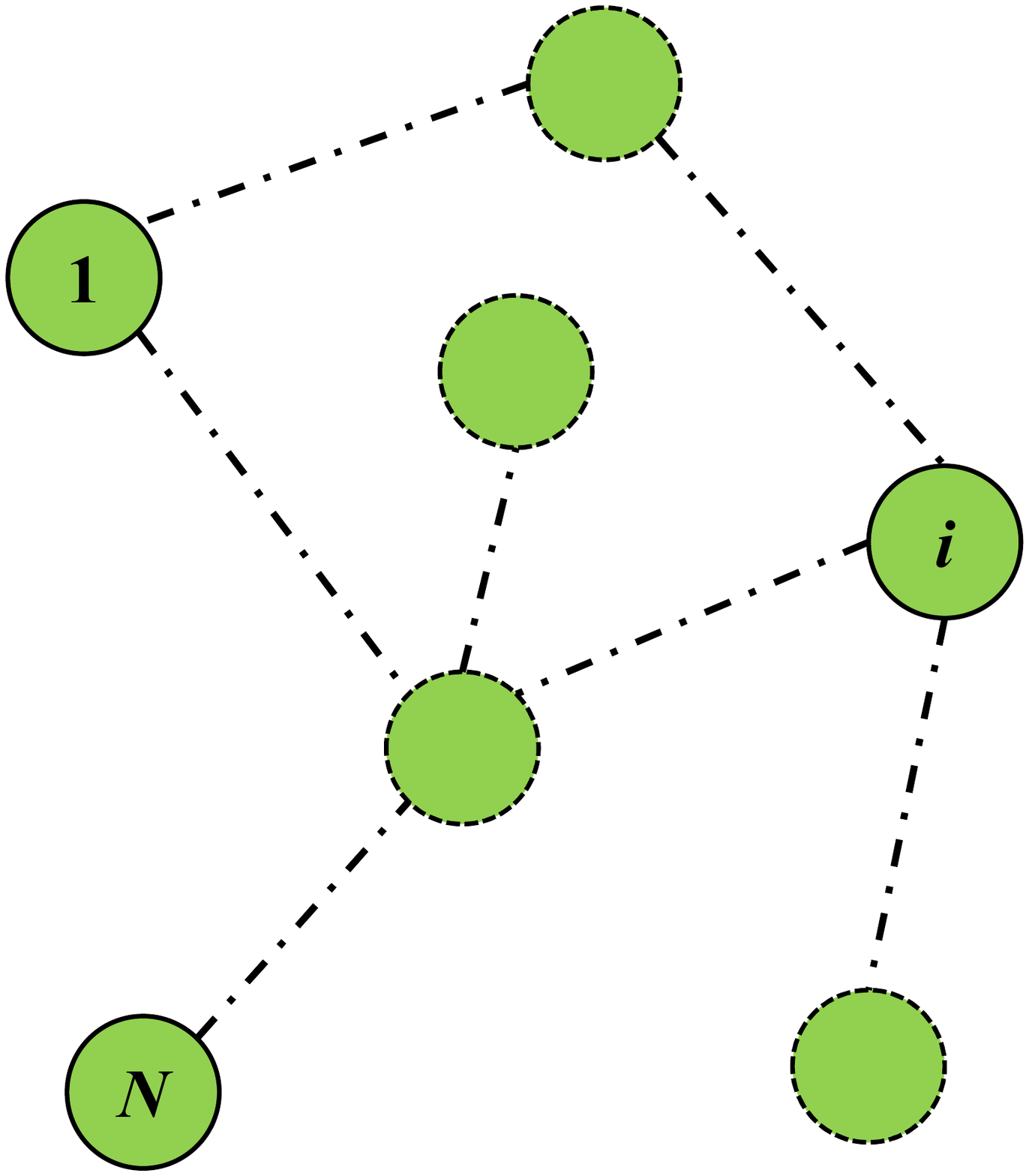}\\
  \caption{An illustration of the communication mechanism of DPG algorithm. Each agent is only linked to the agents who share the same coupling constraints.}\label{an2}
\end{figure}

To realize decentralized computations, we let the updating of ${\bm{\lambda}}_i= [\bm{\theta}_i^{\top}, {\bm{\mu}_i}^{\top}]^{\top} \in \mathds{R}^{B+M}$ be maintained by agent $i$, i.e.,
\begin{align}\label{f33}
 {\bm{\lambda}}_i(k+1) = \mathrm{prox}^{c}_{{q}_i} [{\bm{\lambda}}_i(k) - c \nabla_{{\bm{\lambda}}_i} {P} ({\bm{\lambda}}(k))],
\end{align}
which means
\begin{align}
\bm{\theta}_i(k+1)&  = {\bm{\theta}_i}(k) - c  \nabla_{{\bm{\theta}_i}} {P} ({\bm{\lambda}}(k)), \label {56+1} \\
{\bm{\mu}_i}(k+1)&  = \mathrm{prox}^{c}_{{q}_i} [{\bm{\mu}_i}(k) - c  \nabla_{{\bm{\mu}_i}} {P} ({\bm{\lambda}}(k))]. \label {56+2}
\end{align}
Note that ${\mathbf{F}}_i {\bm{\lambda}} = \bm{\mu}_i$, hence the variables of ${q}_i$ are decoupled from each other. However, each ${p}_i ({\bm{\lambda}})$ contains the information $\sum_{l\in \mathcal{V}} (\mathbf{A}^{(l)}_i )^{\top} \bm{\theta}_l = \sum_{l\in \mathcal{V}_i\cup \{i\}} (\mathbf{A}^{(l)}_i )^{\top} \bm{\theta}_l$ (due to (\ref{37})), which means ${P}({\bm{\lambda}})$ is coupled among the neighbouring agents. Therefore, to compute the complete gradient vector $\nabla_{{\bm{\lambda}}_i} {P} ({\bm{\lambda}}(k))$, agent $i$ needs to collect $\nabla_{{\bm{\lambda}}_i} {p}_l ({\bm{\lambda}}(k))$ from neighbour $l\in \mathcal{V}_i$.
The communication and computation mechanisms of DPG algorithm are shown in Fig. \ref{an2} and Algorithm \ref{ax-1}, respectively.

\begin{algorithm}
\caption{Dual Proximal Gradient Algorithm}\label{ax-1}
\begin{algorithmic}[1]
\State Initialize ${\bm{\lambda}}(0)$ such that ${\Psi}({\bm{\lambda}}(0)) \in \mathds{R}$. Determine step-size $c>0$.
\For {$k= 0,1,2,...$}
\For {$i= 1,2,...,{N}$} (in parallel)
\State Receive $\nabla_{{\bm{\lambda}}_i} {p}_l ({\bm{\lambda}}(k))$ from neighbour $l \in \mathcal{V}_i$.
\State Update ${\bm{\lambda}}_i$ by (\ref{f33}).
\EndFor
\State Obtain an output ${\bm{\lambda}}^{\mathrm{out}}$ under certain convergence criterion.
\EndFor
\end{algorithmic}
\end{algorithm}

\begin{Remark}\label{7}
As seen in Algorithm 1, compared with symmetric scenarios, the asymmetric individual constraints introduce asymmetric Lagrangian multipliers for the coupling constraints, where the dual variables are decomposed in a natural way and no global consensus of $\bm{\theta}_i$ is required.
\end{Remark}

To apply (\ref{f33}), one need to derive (i) $f^{\diamond}_i$ for $\bm{\theta}$ and $\bm{\mu}$ and (ii) the proximal mapping of $(g_i +{\mathds{I}}_{{\Omega}_{i}} )^{\diamond}$ for $\bm{\mu}$, $i \in \mathcal{V}$. For (i), $f^{\diamond}_i$ can be easily obtained if $f_i$ is simple-structured, e.g., $f_i$ is a quadratic function.\cite[Sec. 3.3.1]{boyd2004convex} For (ii), a feasible method is introduced in the following remark, which can avoid the calculation of the proximal mapping of $(g_i +{\mathds{I}}_{{\Omega}_{i}} )^{\diamond}$.

\begin{Remark}\label{rm6}
Based on Lemma \ref{md},
the updating of $\bm{\mu}_i$ in Algorithm \ref{ax-1} can be equivalently written as
\begin{align}
& \bm{\varrho}_i(k)= {\bm{\mu}_i}(k) - c  \nabla_{{\bm{\mu}_i}} {P} ({\bm{\lambda}}(k)), \label{r6} \\
& \bm{\mu}_i(k+1) = \mathrm{prox}^{c}_{{q}_i} [ \bm{\varrho}_i(k) ]  =  \bm{\varrho}_i(k)-  c \mathrm{prox}^{\frac{1}{c}}_{{q}^{\diamond}_i} [\frac{\bm{\varrho}_i(k)}{c} ], \label{r7}
\end{align}
with ${q}^{\diamond}_i({\bm{\lambda}})=(g_i +{\mathds{I}}_{{\Omega}_{i}} )^{\diamond \diamond}({\mathbf{F}}_i{\bm{\lambda}})=(g_i +{\mathds{I}}_{{\Omega}_{i}} )({\mathbf{F}}_i{\bm{\lambda}})$ due to the convexity and lower semi-continuity of $g_i +{\mathds{I}}_{{\Omega}_{i}}$, where $(g_i +{\mathds{I}}_{{\Omega}_{i}} )^{\diamond \diamond}$ is the biconjugate of $g_i +{\mathds{I}}_{{\Omega}_{i}}$.\cite[Sec. 3.3.2]{boyd2004convex} \footnote[4]{(\ref{r6}) can be included in (\ref{r7}) to generate a one-step formula for $\bm{\mu}_i$.} With this arrangement, the calculation of the proximal mapping of $(g_i +{\mathds{I}}_{{\Omega}_{i}})^{\diamond} $ is not required as shown in (\ref{r7}), which reduces the computational complexity when the proximal mapping of $g_i +{\mathds{I}}_{{\Omega}_{i}}$ is easier to obtain by available formulas.\cite[Sec. 6.3]{beck2017first} For example, in some $l_1$ regularization problems (e.g., $g_i(\mathbf{x}_i)= \| \mathbf{x}_i \|_1$, $\Omega_i = \mathds{R}^M$), the proximal mapping of $l_1$-norm is known as the soft thresholding operator with analytical solution.\cite[Sec. 6.3]{beck2017first} In addition, if $g_i=0$ (i.e., smooth cost functions with local constraints), the proximal mapping of ${\mathds{I}}_{{\Omega}_{i}}$ is an Euclidean projection onto ${\Omega}_{i}$.\cite[Sec. 1.2]{parikh2014proximal}
\end{Remark}

Additional to the method in Remark \ref{rm6}, the following remark explains how to implement DPG algorithm for certain general form of $g_i +{\mathds{I}}_{{\Omega}_{i}}$.
\begin{Remark}
If the proximal mapping of $g_i +{\mathds{I}}_{{\Omega}_{i}}$ cannot be obtained efficiently, a feasible method is to construct a strongly convex $g_i$ (e.g., shift a strongly convex component from $f_i$ to $g_i$). By the definition of proximal mapping, (\ref{56+2}) can be rewritten as
\begin{align}\label{f31}
{\bm{\mu}_i}(k+1)  = & \arg \min_{\mathbf{v}} (q_i(\mathbf{v})  + \frac{1}{2c} \| \mathbf{v} - {\bm{\mu}_i}(k) + c  \nabla_{{\bm{\mu}_i}} {P} ({\bm{\lambda}}(k)) \|^2 ).
\end{align}
(\ref{f31}) can be solved with gradient descent method by computing the gradient of $q_i$ with the help of Lemma \ref{l1}, i.e.,
\begin{align}
\nabla_{\mathbf{v}}  q_i(\mathbf{v}) = & \nabla_{\mathbf{v}}(g_i +{\mathds{I}}_{{X}_{i}} )^{\diamond}(\mathbf{v})
=  \arg \max_{\mathbf{u}} (\mathbf{v}^{\top} \mathbf{u} -(g_i + \mathds{I}_{X_i})(\mathbf{u})),
\end{align}
which can be completed with local information. In this case, the DPG algorithm can adapt to general nonsmooth $g_i + \mathds{I}_{X_i}$ with a compromise on an inner-loop optimization process.
\end{Remark}

\subsection{Asynchronous Dual Proximal Gradient Algorithm}\label{as4+2}

In the following, we propose an Asyn-DPG algorithm by extending the usage of DPG algorithm to asynchronous networks.

In synchronous networks, the information accessed by the agents is assumed to be up-to-date, which requires efficient data transmission and can be restrictive for some large-scale networks.\cite{chou1990synchronizing} To address this issue, we propose an Asyn-DPG algorithm for asynchronous networks by considering communication delays.
To this end, based on the setup of Problem (P4), we define $\tau(k)$ as the time instant previous to instant $k$ with $k-\tau(k) \geq 0$.\footnote[5]{In this work, the communication delays are considered to be uniform for the agents at each instant but varying along different instants.\cite{guo2018distributed,wang2013consensus}} Therefore, the accessed dual information at instant $k$ may not be the latest version ${{\bm{\lambda}}}(k)$ but a historical version ${{\bm{\lambda}}}(\tau(k))$. It is reasonable to assume that certain agent always knows the latest information of itself.

\begin{figure}[htbp]
  \centering
  \includegraphics[width=9cm]{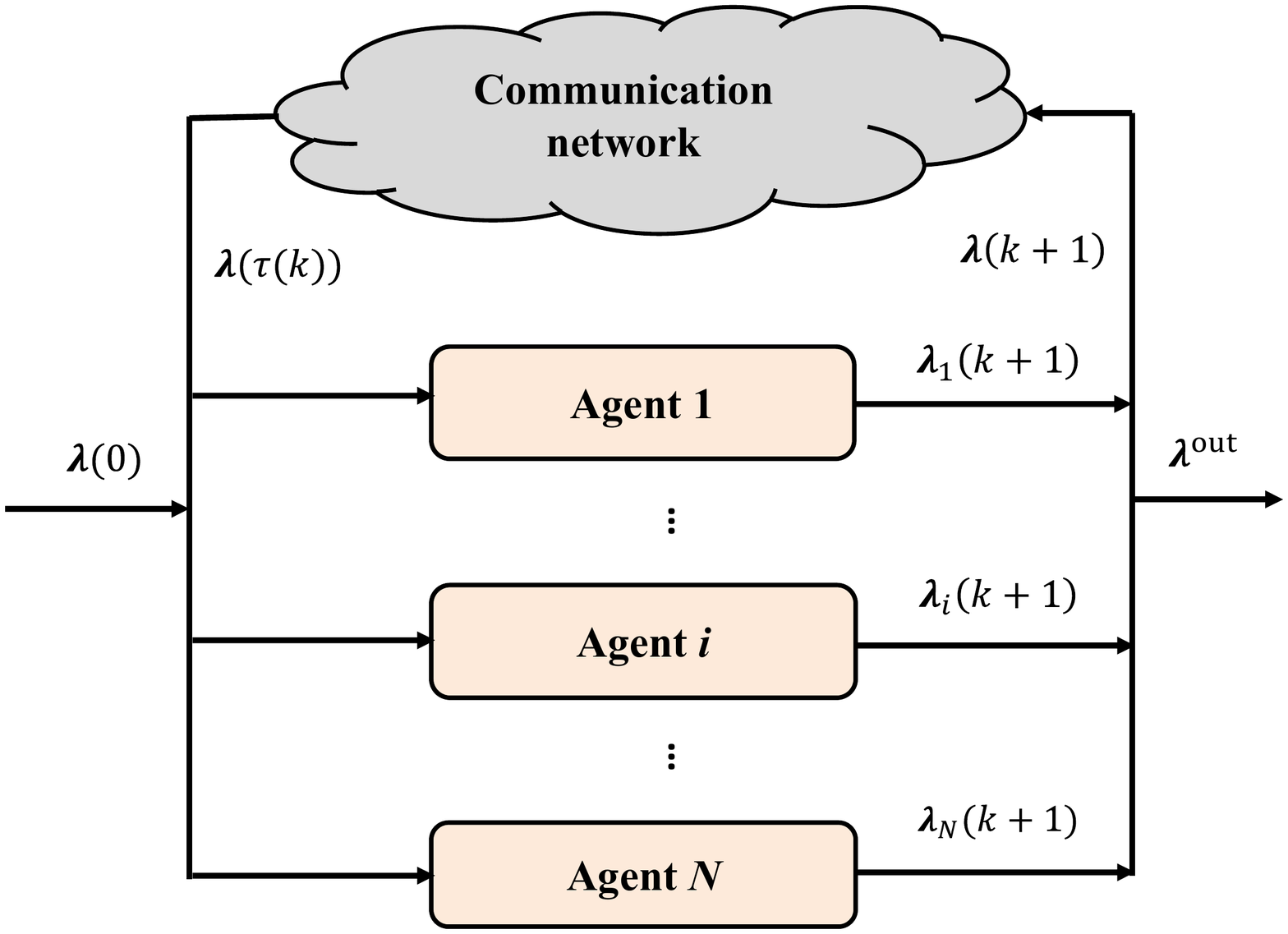}\\
  \caption{An illustration of the computation mechanism of Asyn-DPG algorithm. Each agent pushes the updated state into the network by using possibly delayed information of others.}\label{an1}
\end{figure}

\begin{Assumption}\label{a4}
The communication delays in the network are upper bounded by $D \in \mathds{N}$, which means $0 \leq k-\tau(k) \leq D$, $\tau(k) \in \mathds{N}$, $k \in \mathds{N}$.
\end{Assumption}
The upper bound of delays is a commonly used assumption in asynchronous networks.\cite{wang2013consensus,zhou2018distributed} By allowing for the heterogenous steps-sizes, the proposed Asyn-DPG algorithm is designed as
\begin{align}\label{55}
 {\bm{\lambda}}_i(k+1) = & \mathrm{prox}^{c_i}_{{q}_i} [{\bm{\lambda}}_i(k) -  c_i  \nabla_{{\bm{\lambda}}_i} {P} ({\bm{\lambda}}(\tau(k)))].
\end{align}
The computation mechanism of the Asyn-DPG algorithm is shown in Algorithm \ref{ax1} and Fig. \ref{an1}.

\begin{algorithm}
\caption{Asynchronous Dual Proximal Gradient Algorithm}\label{ax1}
\begin{algorithmic}[1]
\State Initialize ${\bm{\lambda}}(0)$ such that ${\Psi}({\bm{\lambda}}(0)) \in \mathds{R}$. Determine step-size $c_i >0$, $\forall i \in \mathcal{V}$.
\For {$k= 0,1,2,...$}
\For {$i= 1,2,...,{N}$} (in parallel)
\State Receive $\nabla_{{\bm{\lambda}}_i} {p}_l ({\bm{\lambda}} (\tau(k)))$ from neighbour $l\in \mathcal{V}_i$.
\State Update ${\bm{\lambda}}_i$ by (\ref{55}).
\EndFor
\State Obtain an output ${\bm{\lambda}}^{\mathrm{out}}$ under certain convergence criterion.
\EndFor
\end{algorithmic}
\end{algorithm}

Note that (\ref{55}) can be decomposed as indicated in (\ref{56+1}) and (\ref{56+2}), i.e.,
\begin{align}
\bm{\theta}_i(k+1) = & {\bm{\theta}_i}(k) - c_i  \nabla_{{\bm{\theta}_i}} {P} ({\bm{\lambda}}(\tau(k))), \label {} \\
{\bm{\mu}_i}(k+1) = & \mathrm{prox}^{c_i}_{{q}_i} [{\bm{\mu}_i}(k) - c_i  \nabla_{{\bm{\mu}_i}} {P} ({\bm{\lambda}}(\tau(k)))]. \label {}
\end{align}

\begin{Remark}\label{rm7}
Similar to the discussion in Remark \ref{rm6}, the updating of $\bm{\mu}_i$ in Algorithm \ref{ax1} can be rewritten as
\begin{align}
& \bm{\varrho}'_i(k)  =  {\bm{\mu}_i}(k) - c_i \nabla_{{\bm{\mu}_i}} {P} ({\bm{\lambda}}(\tau (k))), \label{r8} \\
& \bm{\mu}_i(k+1)  =  \mathrm{prox}^{c_i}_{{q}_i} [ \bm{\varrho}'_i(k) ] = \bm{\varrho}'_i(k)- c_i \mathrm{prox}^{\frac{1}{c_i}}_{{q}^{\diamond}_i} [\frac{\bm{\varrho}'_i(k)}{c_i} ] \label{r9}
\end{align}
to reduce the computational complexity.
\end{Remark}

\section{Convergence Analysis and Discussion}\label{sa5}

The convergence analysis of the proposed DPG and Asyn-DPG algorithms is conducted in this section.

\begin{Lemma}\label{lam2+1}
With Assumption \ref{a1}, the Lipschitz constant of $\nabla_{{\bm{\lambda}}} {P}({\bm{\lambda}})$ is given by ${h} =  \sum_{i \in \mathcal{V}} {h}_i$, where ${h}_i= \frac{ \| {\mathbf{C}}_i \|^2 }{\sigma_i}$.
\end{Lemma}
See the proof in Appendix \ref{lam2p}.

\begin{Theorem}\label{th1+1}
Suppose that Assumptions \ref{a0}-\ref{as4} hold. Let $\frac{1}{c} = {h}$. By Algorithm \ref{ax-1}, for any ${\bm{\lambda}}^* \in \mathcal{H}$ and $K \in \mathds{N}_+$, we have
\begin{align}\label{}
{\Psi}({\bm{\lambda}} (K)) - {\Psi}({\bm{\lambda}}^*) \leq \frac{{h}\| {\bm{\lambda}} (0)- {\bm{\lambda}}^* \|^2}{2K}.
\end{align}
\end{Theorem}
Note that the structure of (\ref{8}) is consistent with the ISTA algorithm with a constant step-size.\cite{beck2009fast} Therefore, the result of Theorem \ref{th1+1} can be deduced with the existing proof by employing the Lipschitz constant ${h}$.\cite[Thm. 3.1]{beck2009fast} Hence, detailed proof is omitted for simplicity.

\begin{Lemma}\label{lam10}
Based on Assumption \ref{a4}, for certain $K \in \mathds{N}_+$, we have
\begin{align}
\sum^K_{k=0} (\| {\bm{\lambda}} (k+1) - {\bm{\lambda}}(k) \|^2 + \| {\bm{\lambda}} (k) & - {\bm{\lambda}}(k-1) \|^2  + \cdots  + \| {\bm{\lambda}} (\tau(k)+1) - {\bm{\lambda}}(\tau(k)) \|^2) \nonumber \\
& \leq  \sum^K_{k=0}(D+1) \| {\bm{\lambda}} (k+1) - {\bm{\lambda}}(k) \|^2, \label{61} \\
\sum^K_{k=0} k (\| {\bm{\lambda}} (k+1) - {\bm{\lambda}}(k) \|^2 + \| {\bm{\lambda}} (k) & - {\bm{\lambda}}(k-1) \|^2  + \cdots + \| {\bm{\lambda}} (\tau(k)+1) - {\bm{\lambda}}(\tau(k)) \|^2) \nonumber \\
& \leq \sum_{k=0}^K \frac{(2k+D)(D+1)}{2} \| {\bm{\lambda}} (k+1) - {\bm{\lambda}}(k) \|^2. \label{60}
\end{align}
\end{Lemma}
See the proof in Appendix \ref{lam10p}.

\begin{Theorem}\label{ap1}
Suppose that Assumptions \ref{a0}-\ref{a4} hold. By Algorithm \ref{ax1}, given that
\begin{equation}\label{56}
\frac{1}{c_i} \geq {h} (D+1)^2,
\end{equation}
for certain $K \in \mathds{N}_+ \bigcap [ \lceil \frac{D}{2}\rceil, +\infty) $ and any ${\bm{\lambda}}^* \in \mathcal{H}$, we have
\begin{align}\label{57}
{\Psi}({\bm{\lambda}}(K+1)) - {\Psi}({\bm{\lambda}}^*) \leq \frac{\Lambda(c_1,...,c_{N},D)}{K+1},
\end{align}
where $\Lambda(c_1,...,c_{N},D) =  \sum^{\lfloor \frac{D}{2} \rfloor}_{k=0} \sum_{i \in \mathcal{V}} ( \frac{{h} (2k+D)(D+1)^2}{4} - \frac{k}{ c_i } ) \| {\bm{\lambda}}_i(k+1) - {\bm{\lambda}}_i(k) \|^2 + \sum_{i \in \mathcal{V}} \frac{1}{2 c_i} \| {\bm{\lambda}}_i(0)  - {\bm{\lambda}}^*_i \|^2$.
\end{Theorem}
See the proof in Appendix \ref{ap1p}.

\section{Numerical Result}\label{sa6}

In this section, we will verify the feasibility of Algorithms \ref{ax-1} and \ref{ax1} by considering a social welfare optimization problem in an electricity market with 2 utility companies (UCs) and 3 energy users.

\subsection{Simulation Setup}
The social welfare optimization problem of the market is formulated as follows.
\begin{align}\label{}
 \textbf{(P5)} \quad \min_{\mathbf{x}} \quad & \sum_{i \in \mathcal{V}_{\mathrm{UC}}} C_i(x^{\mathrm{UC}}_{i}) - \sum_{j \in \mathcal{V}_{\mathrm{user}}} U_j(x^{\mathrm{user}}_{j}) \nonumber\\
  \hbox{subject to} \quad & \sum_{i \in \mathcal{V}_{\mathrm{UC}}} x^{\mathrm{UC}}_{i} = \sum_{j \in \mathcal{V}_{\mathrm{user}}} x^{\mathrm{user}}_{j}, \label{s1} \\
  & x^{\mathrm{UC}}_{i} \in \Omega^{\textrm{UC}}_i, \quad \forall i \in  \mathcal{V}_{\mathrm{UC}}, \\
  & x^{\mathrm{user}}_{j} \in \Omega^{\textrm{user}}_j, \quad \forall j \in  \mathcal{V}_{\mathrm{user}}.
\end{align}
In Problem (P5), $\mathcal{V}_{\mathrm{UC}}$ and $\mathcal{V}_{\mathrm{user}}$ are the sets of UCs and users, respectively. $\mathbf{x} = [x^{\mathrm{UC}}_{1},...,x^{\mathrm{UC}}_{|\mathcal{V}_{\mathrm{UC}}|}, x^{\mathrm{user}}_{1},...,x^{\mathrm{user}}_{|\mathcal{V}_{\mathrm{user}}|}]^{\top}$ with $x^{\mathrm{UC}}_{i}$ and $x^{\mathrm{user}}_{j}$ being the quantities of energy generation and consumption of UC $i$ and user $j$, respectively. $C_i(x^{\mathrm{UC}}_{i})$ is the cost function of UC $i$ and $U_j(x^{\mathrm{user}}_{j})$ is the utility function of user $j$, $i \in \mathcal{V}_{\mathrm{UC}}$, $j \in \mathcal{V}_{\mathrm{user}}$. The constraint (\ref{s1}) ensures the supply-demand balance in the market. Define constraint matrix $ {\mathbf{A}} =  [\mathbf{1}^{\top}_{|\mathcal{V}_{\mathrm{UC}}|}, -\mathbf{1}^{\top}_{|\mathcal{V}_{\mathrm{user}}|}]$. Then, (\ref{s1}) can be represented by ${\mathbf{A}}\mathbf{x} = 0$. $\Omega^{\textrm{UC}}_i = [0,x^{\mathrm{UC}}_{i,\mathrm{max}}]$ and $\Omega^{\textrm{user}}_j = [0,x^{\mathrm{user}}_{j,\mathrm{max}}]$ are local constraints with $x^{\mathrm{UC}}_{i,\mathrm{max}} >0$ and $x^{\mathrm{user}}_{j,\mathrm{max}}>0$ being the upper bounds of $x^{\mathrm{UC}}_{i}$ and $x^{\mathrm{user}}_{j}$, respectively.

The detailed expressions of $C_i(x^{\mathrm{UC}}_{i})$ and $U_j(x^{\mathrm{user}}_{j})$ are designed as
\begin{align}\label{}
& C_i(x^{\mathrm{UC}}_{i}) = \kappa_i (x^{\mathrm{UC}}_{i})^2 + \vartheta_i x^{\mathrm{UC}}_{i} +\beta_i,  \\
& U_j(x^{\mathrm{user}}_{j}) =
\left\{\begin{array}{ll}
\pi_j x^{\mathrm{user}}_{j} - \varsigma_j (x^{\mathrm{user}}_{j})^2, & x^{\mathrm{user}}_{j} \leq \frac{\pi_j}{2\varsigma_j}, \\
  \frac{\pi_j^2}{4 \varsigma_j}, & x^{\mathrm{user}}_{j} > \frac{\pi_j}{2\varsigma_j},
\end{array}
\right.
\end{align}
where $\kappa_i,\vartheta_i,\beta_i,\pi_j,\varsigma_j$ are parameters, $\forall i \in  \mathcal{V}_{\mathrm{UC}}$, $\forall j \in  \mathcal{V}_{\mathrm{user}}$. The values of the parameters are set in Table I.\cite{pourbabak2017novel}



\begin{table}\label{tm2}
\caption{Parameters of UCs and energy users}
\label{tab2}
\begin{center}
\begin{tabular}{p{3mm}p{8mm}p{7mm}p{3mm}p{8mm}p{7mm}p{9mm}p{8mm}}
\bottomrule
& \multicolumn{3}{l}{$ \quad \quad \quad\quad$ UCs} & \multicolumn{3}{l}{$\quad \quad \quad \quad \quad $ Users} \\
\hline
$i/j$ & $\kappa_i$ & $\vartheta_i$ & $\beta_i$ & $x^{\mathrm{UC}}_{i,\mathrm{max}}$ & $\pi_j$ &  $\varsigma_j$  & $x^{\mathrm{user}}_{j,\mathrm{max}}$ \\
\hline
1 & 0.0031& 8.71 & 0 &150 & 17.17 & 0.0935 & 91.79 \\
2 & 0.0074 & 3.53 & 0 &150 & 12.28 & 0.0417 & 147.29 \\
3 & - & - & - &  - &18.42 & 0.1007 & 91.41 \\
\bottomrule
\end{tabular}
\end{center}
\end{table}

To apply the DPG algorithm, we define $\mathbf{A}^{(i),\mathrm{UC}}$ and $\mathbf{A}^{(j),\mathrm{user}}$
as the asymmetric constraint matrices of UC $i$ and user $j$, respectively. Then, by following the derivation of (\ref{28}), the Lagrangian function of Problem (P5) can be obtained as ${L} (\mathbf{x},\mathbf{z},  {\bm{\theta}},{\bm{\mu}})$,
where ${\bm{\theta}} = [ {\theta}^{\mathrm{UC}}_1,..., {\theta}^{\mathrm{UC}}_{\mid \mathcal{V}_{\mathrm{UC}}\mid}, {\theta}^{\mathrm{user}}_1,..., {\theta}^{\mathrm{user}}_{\mid \mathcal{V}_{\mathrm{user}}\mid}]^{\top}$ and ${\bm{\mu}}=[ \mu^{\mathrm{UC}}_1,..., \mu^{\mathrm{UC}}_{\mid \mathcal{V}_{\mathrm{UC}}\mid},$ ${\mu}^{\mathrm{user}}_1,..., {\mu}^{\mathrm{user}}_{\mid \mathcal{V}_{\mathrm{user}}\mid}]^{\top}$.
See Appendix \ref{sa60} for the detailed expressions of ${L} (\mathbf{x},\mathbf{z}, {\bm{\theta}}, {\bm{\mu}})$, $\mathbf{A}^{(i),\mathrm{UC}}$, and $\mathbf{A}^{(j),\mathrm{user}}$.
With some direct calculations, the optimal solution to Problem (P5) is $\mathbf{x}^*=[0,150,48.5,50.2,51.3]^{\top}$.

\subsection{Simulation Result and Discussion}

\subsubsection{Simulation 1}\label{sm}  To apply Algorithm \ref{ax-1}, we consider a fully connected network since all the agents are involved in supply-demand balance constraint. Due to the different individual interpretations of the global constraint, with some linear transformations introduced in Section \ref{6}, we let $[\mathbf{A}^{(1),\mathrm{UC}},\mathbf{A}^{(2),\mathrm{UC}},\mathbf{A}^{(1),\mathrm{user}},\mathbf{A}^{(2),\mathrm{user}},\mathbf{A}^{(3),\mathrm{user}}] = [\mathbf{T}_1 \mathbf{A},\mathbf{T}_2\mathbf{A}, \mathbf{T}_3 \mathbf{A},\mathbf{T}_4 \mathbf{A},\mathbf{T}_5 \mathbf{A}]$,
where $\mathbf{T}_1=1$, $\mathbf{T}_2=2$, $\mathbf{T}_3=-1$, $\mathbf{T}_4=1$, $\mathbf{T}_5=-1$.

\begin{figure}[htbp]
\centering
\includegraphics[width=8cm]{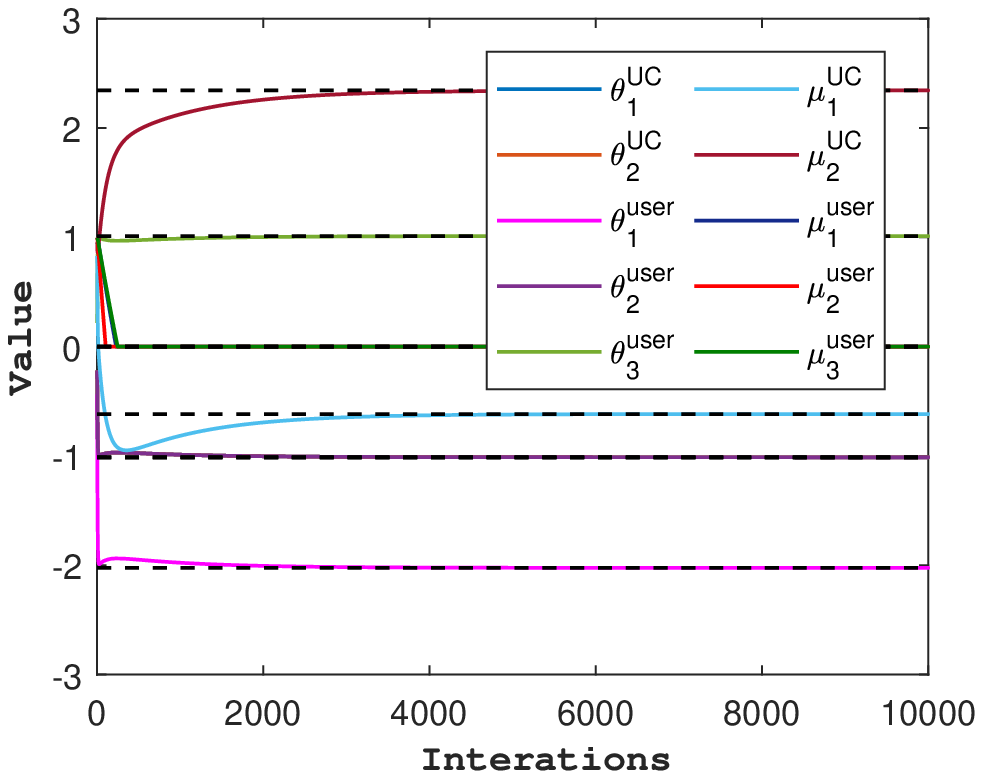}
\caption{Dynamics of ${\bm{\theta}}$ and ${\bm{\mu}}$.}\label{si21}
\end{figure}

\begin{figure}[htbp]
\centering
\includegraphics[width=8cm]{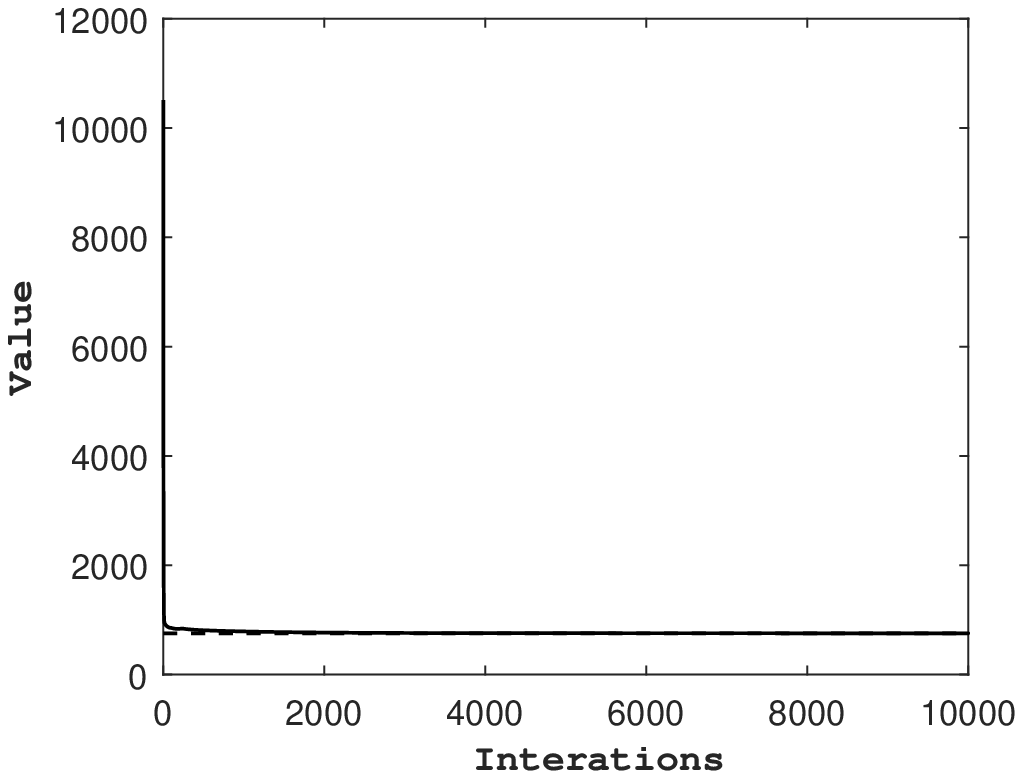}
\caption{Dynamics of ${\Psi}( {\bm{\lambda}})$.}\label{si22}
\end{figure}

The simulation result is shown in Figs. \ref{si21} and \ref{si22}. As shown in Fig. \ref{si21}, $({\bm{\theta}},{\bm{\mu}})$ converges to certain stationary position $({\bm{\theta}}^*,{\bm{\mu}}^*)$. Meanwhile, the optimal solution to the primal problem can be obtained by $\mathbf{x}^*= \arg \min_{\mathbf{x}} {L}( \mathbf{x}, \mathbf{z}, {\bm{\theta}}^*,{\bm{\mu}}^*) = [0,150,48.5,50.2,51.3]^{\top}$, and the value of dual function ${\Psi}( {\bm{\lambda}})$ (as defined in Problem (P4)) converges to the minimum value ${\Psi}^*$, which is around 756.53.

\begin{figure}[htbp]
\centering
\includegraphics[width=8cm]{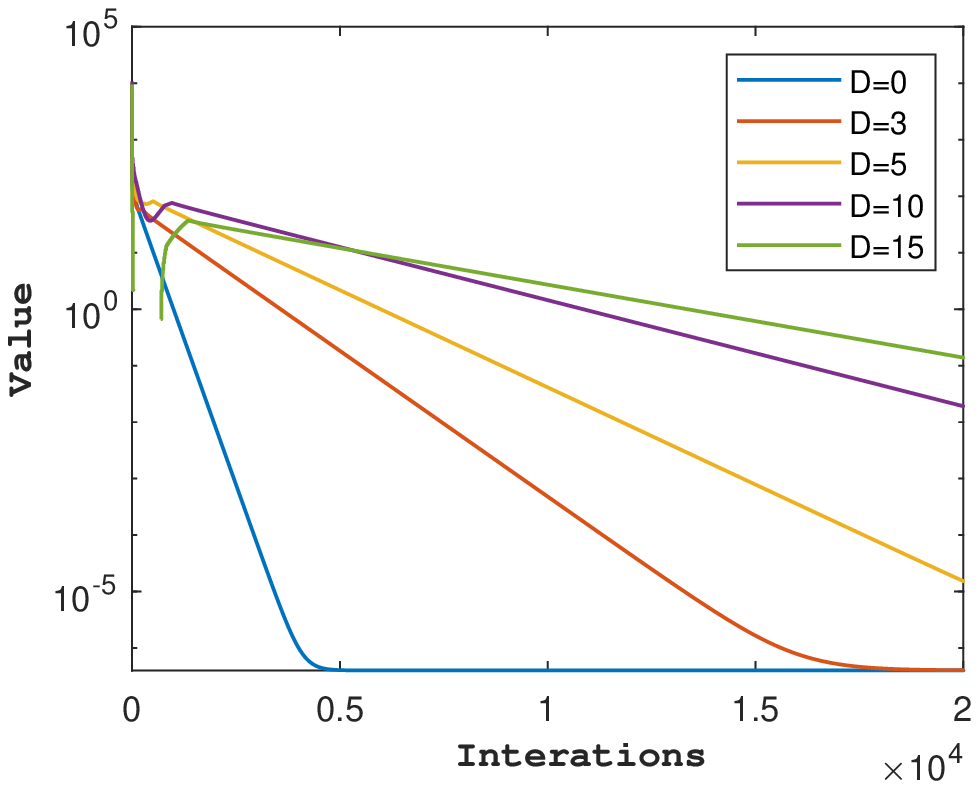}
\caption{Dynamics of $\epsilon$ with different delays.}\label{si23}
\end{figure}

\subsubsection{Simulation 2} To apply Algorithm \ref{ax1}, the upper bound of communication delays is set as $D \in \{0,3,5,10,15\}$. To represent the ``worst delays'', we let $\tau(k)=\max\{0,k-D\}$, $k \in \mathds{N}$. In addition, we define $\epsilon(k)= {\Psi}({\bm{\lambda}}(k)) -{\Psi}^*$ to characterize the dynamics of convergence error.

With the same asymmetric constraints in Simulation 1, the simulation result is shown in Fig. \ref{si23}. It can be seen that, with different delays, the minimum of ${\Psi} ({\bm{\lambda}})$, i.e., ${\Psi}^*$, is achieved asymptotically, which implies the optimal solution to the primal problem is achieved since Simulations 1 and 2 are based on the same setup of Problem (P4). In Fig. \ref{si23}, one can also note that a larger delay can slower the convergence speed, which is consistent with result (\ref{57}), i.e., a larger value of $D$ can produce a larger error bound in certain step. 

\section{Conclusion}\label{sa7}

In this work, we focused on optimizing a class of composite DOPs with both local convex and affine coupling constraints. With different network settings, two dual proximal gradient based algorithms were proposed. As the key feature, all the discussed algorithms resort to the dual problem. Provided that the non-smooth parts of the cost functions are simple-structured, we only need to update dual variables with some simple operations, which leads to the reduction of the overall computational complexity.


%
%
%
%

\appendix{}

\section{Proof of Lemma 3}\label{lam2p}

By the $\frac{1}{\sigma_i}$-Lipschitz continuity of $ \nabla f^{\diamond}_i$ (see Lemma \ref{l1} and Assumption \ref{a1}), we have
\begin{align}\label{z1}
\| \nabla_{\mathbf{u}} f^{\diamond}_i(\mathbf{C}_i\mathbf{u}) - \nabla_{\mathbf{v}} f^{\diamond}_i(\mathbf{C}_i\mathbf{v}) \|
= & \| \mathbf{C}_i^{\top} (\nabla_{\mathbf{C}_i\mathbf{u}} f^{\diamond}_i(\mathbf{C}_i\mathbf{u}) - \nabla_{\mathbf{C}_i\mathbf{v}} f^{\diamond}_i(\mathbf{C}_i\mathbf{v})) \|
\leq  \| \mathbf{C}_i\| \| \nabla_{\mathbf{C}_i\mathbf{u}} f^{\diamond}_i(\mathbf{C}_i\mathbf{u}) - \nabla_{\mathbf{C}_i\mathbf{v}} f^{\diamond}_i(\mathbf{C}_i\mathbf{v}) \|  \nonumber \\
\leq & \frac{\| \mathbf{C}_i\| }{\sigma_i} \| \mathbf{C}_i \mathbf{u} - \mathbf{C}_i \mathbf{v} \|
\leq  \frac{\| \mathbf{C}_i  \|^2 }{\sigma_i}\|  \mathbf{u}- \mathbf{v} \|
=  h_i \|  \mathbf{u}- \mathbf{v} \|,
\end{align}
$\forall \mathbf{u}, \mathbf{v} \in \mathds{R}^{NM+NB}$, which means $ \nabla_{\bm{\lambda}} f^{\diamond}_i(\mathbf{C}_i\bm{\lambda})$ is $h_i$-Lipschitz continuous and, therefore, $\nabla_{\bm{\lambda}} p_i (\bm{\lambda})=\nabla_{\bm{\lambda}} f^{\diamond}(\mathbf{C}_i\bm{\lambda}) + \mathbf{E}_i^{\top}$ is also $h_i$-Lipschitz continuous. Then, the Lipschitz constant of $\nabla_{{\bm{\lambda}}} {P} ({\bm{\lambda}})$ is a linear sum of ${h}_i$, which gives ${h}$, $i \in \mathcal{V}$.

\section{Proof of Lemma 4}\label{lam10p}
For (\ref{61}),
\begin{align}\label{}
& \sum^K_{k=0} (\| {\bm{\lambda}} (k+1) - {\bm{\lambda}}(k) \|^2 + \| {\bm{\lambda}} (k) - {\bm{\lambda}}(k-1) \|^2 + \cdots + \| {\bm{\lambda}} (\tau(k)+1) - {\bm{\lambda}}(\tau(k))\|^2) \nonumber \\
= & (\| {\bm{\lambda}} (K+1) - {\bm{\lambda}}(K) \|^2 + \| {\bm{\lambda}} (K) - {\bm{\lambda}}(K-1) \|^2 + \cdots + \| {\bm{\lambda}} (\tau(K)+1) - {\bm{\lambda}}(\tau(K)) \|^2) \nonumber \\
& + (\| {\bm{\lambda}} (K) - {\bm{\lambda}}(K-1) \|^2   + \| {\bm{\lambda}} (K-1) - {\bm{\lambda}}(K-2) \|^2 + \cdots  + \| {\bm{\lambda}} (\tau(K-1)+1) - {\bm{\lambda}}(\tau(K-1)) \|^2) + \cdots \nonumber \\
& + (\| {\bm{\lambda}} (2) - {\bm{\lambda}}(1) \|^2 + \| {\bm{\lambda}} (1) - {\bm{\lambda}}(0) \|^2)   + \| {\bm{\lambda}} (1) - {\bm{\lambda}}(0) \|^2
\nonumber \\
\leq & \| {\bm{\lambda}} (K+1) - {\bm{\lambda}}(K) \|^2 + 2 \| {\bm{\lambda}} (K) - {\bm{\lambda}}(K-1) \|^2 + \cdots + (D+1) \| {\bm{\lambda}} (\tau(K)+1) - {\bm{\lambda}}(\tau(K)) \|^2 + \cdots  \nonumber \\
&  + (D+1) \| {\bm{\lambda}} (1) - {\bm{\lambda}}(0) \|^2  \nonumber \\
\leq & \sum^K_{k=0}(D+1) \| {\bm{\lambda}} (k+1) - {\bm{\lambda}}(k) \|^2 .
\end{align}
For (\ref{60}),
\begin{align}\label{}
& \sum^K_{k=0} k (\| {\bm{\lambda}} (k+1) - {\bm{\lambda}}(k) \|^2 + \| {\bm{\lambda}} (k) - {\bm{\lambda}}(k-1) \|^2  + \cdots  + \| {\bm{\lambda}} (\tau(k)+1) - {\bm{\lambda}}(\tau(k)) \|^2) \nonumber \\
 = & K( \| {\bm{\lambda}} (K+1) - {\bm{\lambda}}(K) \|^2  + \| {\bm{\lambda}} (K) - {\bm{\lambda}}(K-1) \|^2   + \cdots + \| {\bm{\lambda}} (\tau(K)+1) - {\bm{\lambda}}(\tau(K)) \|^2) \nonumber \\
& + (K-1) (\| {\bm{\lambda}} (K) - {\bm{\lambda}}(K-1) \|^2 + \cdots  + \| {\bm{\lambda}} (\tau(K-1)+1) - {\bm{\lambda}}(\tau(K-1)) \|^2) + \cdots  \nonumber \\
& + 1 \cdot(\| {\bm{\lambda}} (2) - {\bm{\lambda}}(1) \|^2 + \| {\bm{\lambda}} (1) - {\bm{\lambda}}(0) \|^2) + 0 \cdot \| {\bm{\lambda}} (1) - {\bm{\lambda}}(0) \|^2 \nonumber \\
\leq & K \| {\bm{\lambda}} (K+1) - {\bm{\lambda}}(K) \|^2  + ((K-1)+K) \| {\bm{\lambda}} (K) - {\bm{\lambda}}(K-1) \|^2 + \cdots \nonumber \\
& + (k + (k+1) + \cdots + (k+D)) \| {\bm{\lambda}} (k+1) - {\bm{\lambda}}(k) \|^2 + \cdots + (0 + 1 + \cdots + D) \| {\bm{\lambda}} (1) - {\bm{\lambda}}(0) \|^2 \nonumber \\
\leq & \sum_{k=0}^K (k+(k+1) +\cdots+(k+D)) \| {\bm{\lambda}} (k+1) - {\bm{\lambda}}(k) \|^2 \nonumber \\
= & \sum_{k=0}^K \frac{(2k+D)(D+1)}{2} \| {\bm{\lambda}} (k+1) - {\bm{\lambda}}(k) \|^2 .
\end{align}

\section{Proof of Theorem 2}\label{ap1p}
For agent $ i$, by the first-order optimality condition of proximal mapping (\ref{55}), we have
\begin{align}\label{}
\mathbf{0} \in & \nabla_{{\bm{\lambda}}_i} {P}({\bm{\lambda}} (\tau(k))) +  \partial {q}_i({\bm{\lambda}}_i(k+1)) + \frac{1}{c_i }({\bm{\lambda}}_i(k+1) -{\bm{\lambda}}_i(k)).
 \end{align}
By the convexity of ${q}_i$, we have
\begin{align}\label{10+3}
{q}_i  ({\bm{\lambda}}_i  (k+1)) - {q}_i({\bm{\lambda}}_i) \leq & \langle \nabla_{{\bm{\lambda}}_i} {P}({\bm{\lambda}}(\tau(k))),  {\bm{\lambda}}_i - {\bm{\lambda}}_i(k+1) \rangle   + \frac{1}{c_i } \langle {\bm{\lambda}}_i (k+1)- {\bm{\lambda}}_i(k), {\bm{\lambda}}_i - {\bm{\lambda}}_i(k+1) \rangle.
\end{align}
Summing up (\ref{10+3}) over $\mathcal{V}$ gives
\begin{align}\label{10+1}
{Q} ({\bm{\lambda}}(k+1))- {Q}({\bm{\lambda}}) \leq & \langle \nabla_{{\bm{\lambda}}} {P}({\bm{\lambda}}(\tau(k))),  {\bm{\lambda}} - {\bm{\lambda}}(k+1) \rangle + \sum_{i\in \mathcal{V}}\frac{1}{c_i } \langle {\bm{\lambda}}_i (k+1)- {\bm{\lambda}}_i(k), {\bm{\lambda}}_i - {\bm{\lambda}}_i(k+1) \rangle,
\end{align}
where the separability of ${Q}$ is used. By the Lipschitz continuity of $\nabla_{{\bm{\lambda}}} {P}$ and convexity of ${P}$, we have
\begin{align}\label{11}
{P} ({\bm{\lambda}}(k+1)) - {P}({\bm{\lambda}}(\tau(k))) \leq &  \langle \nabla_{{\bm{\lambda}}} {P}({\bm{\lambda}}(\tau(k))),
{\bm{\lambda}} (k+1) - {\bm{\lambda}}(\tau(k))  \rangle + \frac{{h}}{2} \|  {\bm{\lambda}}(k+1) - {\bm{\lambda}}(\tau(k))  \|^2 \nonumber \\
\leq & {P}({\bm{\lambda}})  - \langle \nabla_{{\bm{\lambda}}} {P}({\bm{\lambda}}(\tau(k))),   {\bm{\lambda}} - {\bm{\lambda}} (\tau(k)) \rangle + \langle \nabla_{{\bm{\lambda}}} {P}({\bm{\lambda}}(\tau(k))),   {\bm{\lambda}} (k+1) - {\bm{\lambda}}(\tau(k)) \rangle \nonumber \\
& + \frac{{h}}{2} \| {\bm{\lambda}}(k+1) - {\bm{\lambda}}(\tau(k)) \|^2 \nonumber \\
\leq & {P}({\bm{\lambda}})  + \langle \nabla_{{\bm{\lambda}}} {P}({\bm{\lambda}}(\tau(k))),   {\bm{\lambda}}(k+1) - {\bm{\lambda}} \rangle + \frac{{h}}{2} \| {\bm{\lambda}}(k+1) - {\bm{\lambda}}(\tau(k)) \|^2.
\end{align}
Adding together the both sides of (\ref{10+1}) and (\ref{11}) gives
\begin{align}\label{12}
{\Psi} ({\bm{\lambda}}(k+1)) - {\Psi}({\bm{\lambda}})
\leq & \sum_{i \in \mathcal{V}} \frac{1}{c_i } \langle {\bm{\lambda}}_i(k+1)- {\bm{\lambda}}_i(k),{\bm{\lambda}}_i -{\bm{\lambda}}_i(k+1)  \rangle  + \frac{{h} }{2} \| {\bm{\lambda}}(k+1)
 -{\bm{\lambda}}(\tau(k)) \|^2  \nonumber \\
= & \frac{{h}}{2}  \| {\bm{\lambda}} (k+1) - {\bm{\lambda}}(\tau(k)) \|^2 - \sum_{i \in \mathcal{V}} \frac{1}{2c_i } \| {\bm{\lambda}}_i (k+1) - {\bm{\lambda}}_i(k) \|^2  \nonumber \\
& + \sum_{i \in \mathcal{V}} \frac{1}{2c_i } (\| {\bm{\lambda}}_i(k)  - {\bm{\lambda}}_i \|^2 - \| {\bm{\lambda}}_i(k+1) - {\bm{\lambda}}_i \|^2 ),
\end{align}
where relation $\mathbf{u}^{\top} \mathbf{v} = \frac{1}{2} (\| \mathbf{u} \|^2 + \| \mathbf{v} \|^2 - \| \mathbf{u} - \mathbf{v} \|^2)$ is used, $\forall \mathbf{u},\mathbf{v} \in \mathds{R}^{M+ B}$. By letting ${\bm{\lambda}}={\bm{\lambda}}^*$ in (\ref{12}) and summing up the result over $k=0,...,K$, we have
\begin{align}\label{22}
 \sum_{n=0}^K ({\Psi}({\bm{\lambda}}(k+1)) - {\Psi}({\bm{\lambda}}^*))
\leq & \sum_{k=0}^K ( \frac{{h} }{2}  \| {\bm{\lambda}} (k+1) - {\bm{\lambda}}(\tau(k)) \|^2 -  \sum_{i \in \mathcal{V}} \frac{1}{2c_i } \| {\bm{\lambda}}_i (k+1) -  {\bm{\lambda}}_i(k) \|^2 \nonumber \\
& + \sum_{i \in \mathcal{V}} \frac{1}{2c_i } (\| {\bm{\lambda}}_i(k)  - {\bm{\lambda}}^*_i \|^2 - \| {\bm{\lambda}}_i(k+1)  - {\bm{\lambda}}^*_i \|^2 ) ) \nonumber \\
\leq & \sum_{k=0}^K ( \frac{{h}  (D+1)}{2}  ( \| {\bm{\lambda}} (k+1) - {\bm{\lambda}}(k) \|^2 + \cdots + \| {\bm{\lambda}} (\tau(k)+1) - {\bm{\lambda}}(\tau(k)) \|^2) \nonumber \\
& - \sum_{i \in \mathcal{V}} \frac{1}{2c_i } \| {\bm{\lambda}}_i (k+1)  - {\bm{\lambda}}_i(k) \|^2  +  \sum_{i \in \mathcal{V}} \frac{1}{2c_i } (\| {\bm{\lambda}}_i(k)  - {\bm{\lambda}}^*_i \|^2 - \| {\bm{\lambda}}_i(k+1)  - {\bm{\lambda}}^*_i \|^2 ) ) \nonumber \\
\leq & \sum_{k=0}^K ( \frac{{h} (D+1)^2}{2} \| {\bm{\lambda}} (k+1) - {\bm{\lambda}}(k) \|^2  -  \sum_{i \in \mathcal{V}} \frac{1}{2c_i } \| {\bm{\lambda}}_i (k+1)  - {\bm{\lambda}}_i(k) \|^2 \nonumber \\
& + \sum_{i \in \mathcal{V}} \frac{1}{2c_i } (\| {\bm{\lambda}}_i(k)  - {\bm{\lambda}}^*_i \|^2 - \| {\bm{\lambda}}_i(k+1)  - {\bm{\lambda}}^*_i \|^2 ) ) \nonumber \\
\leq & \sum_{k=0}^K \sum_{i \in \mathcal{V}} ( \frac{{h}  (D+1)^2}{2} - \frac{1}{2c_i }) \| {\bm{\lambda}}_i(k+1) - {\bm{\lambda}}_i(k) \|^2 +  \sum_{i \in \mathcal{V}} \frac{1}{2c_i } \| {\bm{\lambda}}_i(0)  - {\bm{\lambda}}^*_i \|^2,
\end{align}
where Cauchy-Schwarz inequality and (\ref{61}) are used in the second and third inequalities, respectively. Letting ${\bm{\lambda}}={\bm{\lambda}}(k)$ in (\ref{12}) gives
\begin{align}\label{23}
{\Psi} ({\bm{\lambda}}(k+1)) - {\Psi}({\bm{\lambda}}(k))
\leq & \frac{{h} }{2}  \| {\bm{\lambda}} (k+1) - {\bm{\lambda}}(\tau(k)) \|^2 -  \sum_{i \in \mathcal{V}}  \frac{1}{c_i} \| {\bm{\lambda}}_i(k+1)  - {\bm{\lambda}}_i(k) \|^2.
\end{align}
Multiplying (\ref{23}) by $k$ and summing up the result over $k=0,1,...,K$ gives
\begin{align}\label{24}
\sum_{k=0}^K k({\Psi}({\bm{\lambda}}(k+1)) - {\Psi}({\bm{\lambda}}(k)))
= & \sum_{k=0}^K ( (k+1){\Psi}({\bm{\lambda}}(k+1)) - k{\Psi}({\bm{\lambda}}(k)) - {\Psi}({\bm{\lambda}}(k+1))) \nonumber \\
= & (K+1){\Psi}({\bm{\lambda}}(K+1)) - \sum_{n=0}^K {\Psi}({\bm{\lambda}}(k+1)) \nonumber \\
\leq & \sum^K_{k=0} k ( \frac{{h} }{2} \| {\bm{\lambda}}(k+1) - {\bm{\lambda}}(\tau(k)) \|^2 - \sum_{i \in \mathcal{V}} \frac{1}{c_i} \| {\bm{\lambda}}_i(k+1)  - {\bm{\lambda}}_i(k) \|^2 ) \nonumber \\
\leq & \sum^K_{k=0} k ( \frac{{h}  (D+1)}{2}  ( \| {\bm{\lambda}} (k+1) - {\bm{\lambda}}(k) \|^2 + \cdots + \| {\bm{\lambda}} (\tau(k)+1) - {\bm{\lambda}}(\tau(k)) \|^2)  \nonumber \\
& - \sum_{i \in \mathcal{V}} \frac{1}{c_i} \| {\bm{\lambda}}_i(k+1)  - {\bm{\lambda}}_i(k) \|^2 ) \nonumber \\
\leq & \sum^K_{k=0} \sum_{i \in \mathcal{V}}  ( \frac{{h} (2k+D)(D+1)^2}{4} - \frac{k}{ c_i } )  \| {\bm{\lambda}}_i(k+1) - {\bm{\lambda}}_i(k) \|^2,
\end{align}
where (\ref{60}) is used in the last inequality. By adding the both sides of (\ref{22}) and (\ref{24}) together, we have
\begin{align}\label{111}
( K+1)({\Psi}({\bm{\lambda}}(K+1)) - {\Psi}({\bm{\lambda}}^*))
 \leq &  \sum_{k=0}^K \sum_{i \in \mathcal{V}} ( \frac{{h}  (D+1)^2}{2} - \frac{1}{2c_i }) \| {\bm{\lambda}}_i(k+1) - {\bm{\lambda}}_i(k) \|^2 \nonumber \\
 & + \sum^K_{k=0} \sum_{i \in \mathcal{V}} ( \frac{{h} (2k+D)(D+1)^2}{4} - \frac{k}{ c_i } )  \| {\bm{\lambda}}_i(k+1) - {\bm{\lambda}}_i(k) \|^2 \nonumber + \sum_{i \in \mathcal{V}} \frac{1}{2 c_i} \| {\bm{\lambda}}_i(0)  - {\bm{\lambda}}^*_i \|^2 \nonumber \\
= & \sum_{k=0}^K \sum_{i \in \mathcal{V}} ( \frac{{h}  (D+1)^2}{2} - \frac{1}{2c_i }) \| {\bm{\lambda}}_i(k+1) - {\bm{\lambda}}_i(k) \|^2 \nonumber \\
& + \sum^{\lfloor \frac{D}{2} \rfloor}_{k=0} \sum_{i \in \mathcal{V}} ( \frac{{h} (2k+D)(D+1)^2}{4} - \frac{k}{ c_i } ) \| {\bm{\lambda}}_i(k+1) - {\bm{\lambda}}_i(k) \|^2  \nonumber \\
 & +  \sum^{K}_{k=\lceil \frac{D}{2} \rceil} \sum_{i \in \mathcal{V}} ( \frac{{h} (2k+D)(D+1)^2}{4} - \frac{k}{ c_i } )  \| {\bm{\lambda}}_i(k+1) - {\bm{\lambda}}_i(k) \|^2 + \sum_{i \in \mathcal{V}} \frac{1}{2 c_i} \| {\bm{\lambda}}_i(0)  - {\bm{\lambda}}^*_i \|^2 \nonumber \\
= & \sum_{k=0}^K \sum_{i \in \mathcal{V}} ( \underbrace{\frac{{h}  (D+1)^2}{2} - \frac{1}{2c_i }}_{= \varpi_1}) \| {\bm{\lambda}}_i(k+1) - {\bm{\lambda}}_i(k) \|^2 \nonumber \\
& +\sum^{\lfloor \frac{D}{2} \rfloor}_{k=0} \sum_{i \in \mathcal{V}} ( \frac{{h} (2k+D)(D+1)^2}{4} - \frac{k}{ c_i } ) \| {\bm{\lambda}}_i(k+1) - {\bm{\lambda}}_i(k) \|^2  + \sum^{K}_{k=\lceil \frac{D}{2} \rceil} \sum_{i \in \mathcal{V}} ( k(\underbrace{\frac{{h}(D+1)^2}{2} - \frac{1}{ 2c_i} }_{= \varpi_2}) \nonumber \\
&  +  (\underbrace{\frac{{h}D(D+1)^2}{4} - \frac{k}{2c_i}}_{= \varpi_3}) ) \| {\bm{\lambda}}_i(k+1) - {\bm{\lambda}}_i(k) \|^2  + \sum_{i \in \mathcal{V}} \frac{1}{2 c_i} \| {\bm{\lambda}}_i(0)  - {\bm{\lambda}}^*_i \|^2 \nonumber \\
\leq & \sum^{\lfloor \frac{D}{2} \rfloor}_{k=0} \sum_{i \in \mathcal{V}} ( \frac{{h} (2k+D)(D+1)^2}{4} - \frac{k}{ c_i } ) \| {\bm{\lambda}}_i(k+1) - {\bm{\lambda}}_i(k) \|^2 + \sum_{i \in \mathcal{V}} \frac{1}{2 c_i} \| {\bm{\lambda}}_i(0)  - {\bm{\lambda}}^*_i \|^2 \nonumber \\
= & \Lambda(c_1,...,c_{N},D),
\end{align}
where $\varpi_1,\varpi_2,\varpi_3\leq 0$ with $\frac{1}{c_i} \geq {h} (D+1)^2$, $i \in \mathcal{V}$. This proves (\ref{57}) .

\section{Matrices and Lagrangian Function in Section 6}\label{sa60}

The asymmetric constraint matrices of UC $i$ and user $j$ are given by
\begin{align}\label{}
\mathbf{A}^{(i),\mathrm{UC}}= & [A_1^{(i),\mathrm{UC}},...,A_{\mid \mathcal{V}_{\mathrm{UC}}\mid }^{(i),\mathrm{UC}},  A_{\mid \mathcal{V}_{\mathrm{UC}}\mid + 1 }^{(i),\mathrm{UC}},... ,A_{\mid \mathcal{V}_{\mathrm{UC}}\mid + \mid \mathcal{V}_{\mathrm{user}}\mid}^{(i),\mathrm{UC}}], \label{aa1}\\
\mathbf{A}^{(j),\mathrm{user}}= & [A_1^{(j),\mathrm{user}},...,A_{\mid \mathcal{V}_{\mathrm{UC}}\mid }^{(j),\mathrm{user}}, A_{\mid \mathcal{V}_{\mathrm{UC}}\mid + 1 }^{(j),\mathrm{user}},... ,A_{\mid \mathcal{V}_{\mathrm{UC}}\mid + \mid \mathcal{V}_{\mathrm{user}}\mid}^{(j),\mathrm{user}}]. \label{aa2}
\end{align}
Similar to the derivation of (\ref{28}), based on (\ref{aa1}) and (\ref{aa2}), one can have
\begin{align}\label{}
 {L} (\mathbf{x},  \mathbf{z}, {\bm{\theta}},{\bm{\mu}}) = &
\sum_{i \in \mathcal{V}_{\mathrm{UC}}} (C_i(x^{\mathrm{UC}}_{i}) + \mathds{I}_{\Omega^{\textrm{UC}}_i} (z^{\mathrm{UC}}_{i}))  + \sum_{j \in \mathcal{V}_{\mathrm{user}}} ( -U_j(x^{\mathrm{user}}_{j})+ \mathds{I}_{\Omega^{\textrm{user}}_j} (z^{\mathrm{user}}_{j}))  \nonumber \\
 &  + \sum_{i \in \mathcal{V}_{\mathrm{UC}}}  x^{\mathrm{UC}}_i ( \sum_{l\in \mathcal{V}_{\mathrm{UC}} } {A}^{(l),\mathrm{UC}}_i {\theta}^{\mathrm{UC}}_l +  \sum_{l'\in \mathcal{V}_{\mathrm{user}} } {A}^{(l'),\mathrm{user}}_i {\theta}^{\mathrm{user}}_{l'}  + {\mu}^{\mathrm{UC}}_i ) \nonumber \\
 & + \sum_{j \in \mathcal{V}_{\mathrm{user}}}  x^{\mathrm{user}}_j ( \sum_{l\in \mathcal{V}_{\mathrm{UC}} } {A}^{(l),\mathrm{UC}}_{\mid \mathcal{V}_{\mathrm{UC}} \mid+j} {\theta}^{\mathrm{UC}}_l +  \sum_{l'\in \mathcal{V}_{\mathrm{user}} } {A}^{(l'),\mathrm{user}}_{\mid \mathcal{V}_{\mathrm{UC}} \mid+j} {\theta}^{\mathrm{user}}_{l'}  + {\mu}^{\mathrm{user}}_j ) \nonumber \\
 &  - \sum_{i \in \mathcal{V}_{\mathrm{UC}}} z^{\mathrm{UC}}_i {\mu}^{\mathrm{UC}}_i - \sum_{j \in \mathcal{V}_{\mathrm{user}}} z^{\mathrm{user}}_j {\mu}^{\mathrm{user}}_j.
\end{align}

%
%
%








\bibliography{1myref}%



\end{document}